\documentclass[dvips,useAMS,printer]{gSSR2e}

\def\P{\mathbb{P}}
\def\R{\mathbb{R}}

\def\vt{\vartheta}
\def\ve{\varepsilon}
\def\vr{\varrho}
\def\f{\varphi}
\def\hvt{\hat\vt}
\def\hf{\hat\f}
\def\cQ{\mathcal{Q}}
\def\cR{\mathcal{R}}
\def\dc{\dot\chi}
\def\dr{\dot\vr}
\def\ds{\dot\sigma}
\def\sj{\sum_{j=1}^N}
\def\avj{\frac{1}{N}\sj}
\def\op{o_p(1)}
\def\opn{o_p(n^{-1/2})}

\newcommand\bel[1]{\begin{equation}\label{#1}}
\def\ee{\end{equation}}

\def\ben{\begin{eqnarray}}
\def\non{\end{eqnarray}}

\def\be*{\begin{eqnarray*}}
\def\no*{\end{eqnarray*}}

\newcommand{\Ref}[1]{(\ref{#1})}

\makeatletter
\def\ps@plain{%
     \let\@mkboth\@gobbletwo
     \let\@oddfoot\@evenfoot
     \def\@oddhead{\hbox to \textwidth{{\small\begin{tabular}[t]{@{}l@{}}
                                                                                     \\ 
                          \end{tabular}}%
                          \hfill%
                          }}%
 \def\@oddfoot{\vbox{}}
     \let\@evenhead\@oddhead
}

\usepackage{hypernat}
   \def\MR#1{\href{http://www.ams.org/mathscinet-getitem?mr=#1}{MR#1}}

\RequirePackage{comment}
\RequirePackage{hyperref}
\excludecomment{frontmatter}

\begin{document}

\begin{frontmatter}

\title{Optimality of estimators for misspecified semi-Markov models}

\begin{aug}
\author[A]{\fnms{Ursula U.} \snm{M\"uller}}
\author[B]{\fnms{Anton} \snm{Schick}}
\author[C]{\fnms{Wolfgang} \snm{Wefelmeyer}\corref{}\ead[label=e1]{wefelm@math.uni-koeln.de}}
\address[A]{Department of Statistics,
  Texas A\&M University,
  College Station, TX 77843-3143, USA}
\address[B]{Department of Mathematical Sciences,
  Binghamton University,
  Binghamton, NY 13902-6000, USA}
\address[C]{Mathematisches Institut,
  Universit\"at zu K\"oln,
  Weyertal 86-90,
  50931 K\"oln, Germany}
\affiliation{Texas A\&M University, Binghamton University and Universit\"at zu K\"oln}
\end{aug}

\begin{abstract}
Suppose we observe a geometrically ergodic semi-Markov process
and have a parametric model for the transition distribution
of the embedded Markov chain,
for the conditional distribution of the inter-arrival times, or for both.
The first two models for the process are semiparametric,
and the parameters can be estimated by conditional maximum likelihood
estimators. The third model for the process is parametric,
and the parameter can be estimated by an unconditional
maximum likelihood estimator.
We determine heuristically the asymptotic distributions of these estimators
and show that they are asymptotically efficient.
If the parametric models are not correct,
the (conditional) maximum likelihood estimators estimate the parameter
that maximizes the Kullback--Leibler information.
We show that they remain asymptotically efficient in a nonparametric sense.
\end{abstract}



\end{frontmatter}

\doi{10.1080/1744250YYxxxxxxxx}
 \issn{1744-2516}
\issnp{1744-2508} \jvol{00} \jnum{00} \jyear{2006} \jmonth{December}

\markboth{Ursula U. M\"uller, Anton Schick and Wolfgang Wefelmeyer}{
Optimality of estimators for misspecified semi-Markov models}

\title{Optimality of estimators for misspecified semi-Markov models}

\author{URSULA U. M\"ULLER$\dag$,
ANTON SCHICK$^1$$\ddag$
\thanks{$^1$Supported in part by NSF Grant DMS0405791}
and WOLFGANG WEFELMEYER$^\ast$$\S$
\thanks{$^\ast$Corresponding author. Email: wefelm@math.uni-koeln.de
\vspace{6pt}} \\
\vspace{6pt}
$\dag$ Department of Statistics,
  Texas A\&M University,
  College Station, TX 77843-3143, USA \\
$\ddag$ Department of Mathematical Sciences,
  Binghamton University,
  Binghamton, NY 13902-6000, USA \\
$\S$ Mathematisches Institut,
  Universit\"at zu K\"oln,
  Weyertal 86-90,
  50931 K\"oln, Germany \\
\vspace{6pt}
\received{v3.1 released December 2006} }

\maketitle

\begin{abstract}
Suppose we observe a geometrically ergodic semi-Markov process
and have a parametric model for the transition distribution
of the embedded Markov chain,
for the conditional distribution of the inter-arrival times, or for both.
The first two models for the process are semiparametric,
and the parameters can be estimated by conditional maximum likelihood
estimators. The third model for the process is parametric,
and the parameter can be estimated by an unconditional
maximum likelihood estimator.
We determine heuristically the asymptotic distributions of these estimators
and show that they are asymptotically efficient.
If the parametric models are not correct,
the (conditional) maximum likelihood estimators estimate the parameter
that maximizes the Kullback--Leibler information.
We show that they remain asymptotically efficient in a nonparametric sense.
\bigskip

\bigskip

\begin{classcode}
Primary: 62M09; secondary: 62F12, 62G20.
\end{classcode}
\end{abstract}

\section{Introduction}

For i.i.d.\ observations, Daniels \cite{Da61} and Huber \cite{Hu67}
show that the maximum likelihood estimator of a misspecified parametric
model estimates the parameter that maximizes the Kullback--Leibler
(\emph{KL}) information, and determine its asymptotic distribution.
Weaker conditions are given by Pollard \cite{Po85}. For applications
see also White \cite{Wh82}, M\"uller \cite{Mu07}, and Doksum,
Ozeki, Kim and Neto \cite{DOKN}.
Analogous results are obtained for parametric Markov chain models
by Ogata \cite{Og80}, for parametric time series by Hosoya \cite{Ho89}
and by Andrews and Pollard \cite{AP94}, and for parametric diffusion models
by McKeague \cite{Mc84} and Kutoyants \cite{Ku88}.
We refer also to the monograph of Kutoyants \cite{Ku04}. Applications
to time series models in econometrics are studied by White \cite{Wh84}
and Sin and White \cite{SW96}, and in the monograph of White \cite{Wh94}.

Greenwood and Wefelmeyer \cite{GW97} prove that the maximum likelihood
estimator of a misspecified parametric Markov chain model is efficient
in a nonparametric sense. Related efficiency results for misspecified
parametric time series are in Dahlhaus and Wefelmeyer \cite{DW96}.
Here we outline corresponding results for semi-Markov processes.
We consider both parametric and semiparametric misspecified models.
The arguments are heuristic; sufficient regularity conditions can be obtained
as in the above references.

Suppose we observe a semi-Markov process $Z_t$, $t\geq 0$,
with values in an arbitrary measurable space $E$, on a time interval
$0\leq t\leq n$. Let $(X_0,T_0),(X_1,T_1),\dots$ denote the embedded
Markov renewal process. Its transition distribution factors as
\[
S(x,dy,du) = Q\otimes R(x,dy,du) = Q(x,dy)R(x,y,du),
\]
where $Q(x,dy)$ is the transition distribution of the embedded Markov chain
$X_0,X_1,\dots$, and $R(x,y,du)$ is the conditional distribution
of the inter-arrival time $U_j=T_j-T_{j-1}$ given $X_{j-1}=x$ and $X_j=y$.

We assume that the embedded Markov chain is stationary.
We write $P_1(dx)$, $P_2(dx,dy)$ and $P_3(dx,dy,du)$ for the stationary laws
of $X_{j-1}$, $(X_{j-1},X_j)$ and $(X_{j-1},X_j,U_j)$, respectively.
Of course,
$P_2 = P_1\otimes Q$ and $P_3 = P_2\otimes R = P_1\otimes Q\otimes R$.
Set $N=\max\{j:T_j\leq n\}$. We note that studying a semi-Markov process
is equivalent to studying the embedded Markov renewal process.
The latter is a Markov chain. Observing the semi-Markov process
up to time $n$ is equivalent to observing the embedded Markov renewal process
up to the random time $N$.

Natural estimators for $P_1$, $P_2$ and $P_3$
are the empirical distributions
\[
\P_1 = \avj\delta_{X_{j-1}}, \qquad
\P_2 = \avj\delta_{(X_{j-1},X_j)}, \qquad
\P_3 = \avj\delta_{(X_{j-1},X_j,U_j)}.
\]
where $\delta_x$ denotes the Dirac measure at a point $x$.

Let $\Theta$ be an open subset of $\R^d$. We consider the following
three models for the semi-Markov process. In \emph{Model Q}
we assume a parametric form $Q=Q_\vt$, $\vt\in\Theta$,
of the transition distribution of the embedded Markov chain. These models
are also considered in Greenwood, M\"uller and Wefelmeyer \cite{GMW04}.
In \emph{Model R} we assume a parametric form $R=R_\vt$, $\vt\in\Theta$,
of the conditional distribution of the inter-arrival times.
In \emph{Model S} we assume parametric forms $Q=Q_\vt$ and $R=R_\vt$,
$\vt\in\Theta$, for both. Of course, the last model covers the case
that $Q$ and $R$ carry different parameters.
We assume that $Q_\vt(x,dy)$ has a density $q_\vt(x,y)$ with respect
to some dominating measure $\mu(dy)$, and $R_\vt(x,y,du)$
has a density $r_\vt(x,y,u)$ with respect
to some dominating measure $\nu(du)$.

If Model Q holds, then the transition distribution of the semi-Markov process
is semiparametric, $S=Q_\vt\times R$, with $R$ an infinite-dimensional
nuisance parameter. A natural estimator of $\vt$ is the
\emph{partial maximum likelihood estimator} $\hvt_Q$, which maximizes
\[
\P_2[\log q_\vt] = \avj\log q_\vt(X_{j-1},X_j).
\]
Suppose that Model Q is misspecified, and that the true transition
distribution of the embedded Markov chain is $Q$.
Then $\P_2[\log q_\vt]$ is an empirical version of the \emph{KL information}
$P_2[\log q_\vt]$. Let $K_Q(P_2)$ denote the parameter that maximizes
$P_2[\log q_\vt]$. We call $K_Q$ a \emph{KL functional}.
Note that the partial maximum likelihood estimator is the empirical version
of the KL functional, $\hvt_Q=K_Q(\P_2)$.
Since Model Q is misspecified, the semi-Markov model is nonparametric.
The empirical distribution $\P_2$ is efficient for $P_2$ in a certain sense.
If the KL functional is smooth, i.e.\ compactly differentiable
in an appropriate sense, it follows that $\hvt_Q=K_Q(\P_2)$
is efficient for $K_Q(P_2)$. We will not use this approach in this paper.
Instead we derive, in Section 3, a stochastic expansion of $\hvt_Q$,
and determine its influence function. We also show that
the KL functional $K_Q$ is pathwise differentiable, and determine
its canonical gradient.
To keep the exposition simple, we do not give regularity conditions
for these results. They can be adapted e.g.\ from those
of Greenwood and Wefelmeyer \cite{GW97}.
It turns out that the  canonical gradient equals the influence function
of $\hvt_Q$. By the characterisation of efficient estimators in Section 2,
this shows that $\hvt_Q$ is efficient in the nonparametric
semi-Markov model. We also show that $\hvt_Q$ remains efficient
when Model Q is true. The advantage of our approach is that we do not need
to check compact differentiability of $K_Q$
and a corresponding efficiency property of $\P_2$.

The other two models are treated analogously. If Model R holds,
then the transition distribution of the semi-Markov process is semiparametric,
$S=Q\otimes R_\vt$, with $Q$ an infinite-dimensional nuisance parameter.
A natural estimator of $\vt$ is the
\emph{partial maximum likelihood estimator} $\hvt_R$, which maximizes
\[
\P_3[\log r_\vt] = \avj\log r_\vt(X_{j-1},X_j,U_j).
\]
Suppose that Model Q is misspecified, and that the true
conditional distribution of the inter-arrival times is $R$.
Then $\P_3[\log r_\vt]$ is an empirical version of $P_3[\log r_\vt]$.
Again we call the latter \emph{KL information}.
We denote by $K_R(P_3)$ the parameter that maximizes $P_3[\log r_\vt]$,
and we call $K_R$ a \emph{KL functional}. Then $\hvt_R=K_R(\P_3)$.
In Section 4 we derive a stochastic expansion of $\hvt_R$
and the canonical gradient of $K_R$ and show that $\hvt_R$ is efficient
in the nonparametric semi-Markov model. We also show that $\hvt_R$
remains efficient when Model R is true.

If Model S holds, then the transition distribution
of the semi-Markov process is parametric, $S_\vt=Q_\vt\otimes R_\vt$.
Set
\[
s_\vt(x,y,u) = q_\vt(x,y)r_\vt(x,y,u).
\]
A natural estimator of $\vt$ is the
\emph{maximum likelihood estimator} $\hvt_S$, which maximizes
\[
\P_3[\log s_\vt] = \P_2[\log q_\vt] + \P_3[\log r_\vt]
= \avj\log q_\vt(X_{j-1},X_j) + \avj\log r_\vt(X_{j-1},X_j,U_j).
\]
Suppose that Model Q is misspecified, and that the true
transition distribution of the embedded Markov renewal process is
$S=Q\otimes R$. Then $\P_3[\log s_\vt]$ is an empirical version
of $P_3[\log s_\vt]$. Again we call the latter \emph{KL information}.
We denote by $K_S(P_3)$ the parameter that maximizes $P_3[\log s_\vt]$,
and we call $K_S$ a \emph{KL functional}. Then $\hvt_S=K_S(\P_3)$.
In Section 5 we derive a stochastic expansion of $\hvt_S$
and the canonical gradient of $K_S$ and show that $\hvt_S$ is efficient
in the nonparametric semi-Markov model. We also show that $\hvt_S$
remains efficient when Model S is true. Section 6 contains some additional
comments.

\section{Characterization of efficient estimators}

We assume that the embedded Markov chain is positive Harris recurrent
and geometrically ergodic in $L_2(P_2)$. We make the usual assumption
that the conditional distribution of the inter-arrival times
does not charge zero.
We also assume that the mean inter-arrival time $m=EU_j$ is finite. Then
\bel{N}
n/N \to m \quad a.s.
\ee
For a function $f\in L_2(P_3)$ we have the strong law of large numbers
\bel{lln}
\avj f(X_{j-1},X_j,U_j) \to P_3[f] \quad \mbox{a.s.}
\ee
For a function $f\in L_2(P_3)$ with $Sf=0$ we have the martingale central limit
theorem
\bel{clt}
n^{-1/2}\sj f(X_{j-1},X_j,U_j) \Rightarrow m^{-1/2}(P_3[f^2])^{1/2} Y,
\ee
where $Y$ denotes a standard normal random variable.

In order to characterize efficient estimators for functionals
of semi-Markov models, we consider a family $Q_\delta$, $\delta\in\Delta$,
of transition distributions of the embedded Markov chain,
and a family $R_\delta$, $\delta\in\Delta$, of conditional distributions
of the inter-arrival time. Here $\Delta$ is a possibly infinite-dimensional
set, the \emph{parameter space}.
We fix $\delta\in\Delta$ and set $Q=Q_\delta$, $R=R_\delta$ and
\[
V = \{v\in L_2(P_2):Qv=0\}, \qquad
W  = \{w\in L_2(P_3):Rw=0\}.
\]
Note that $V$ and $W$ can be viewed as orthogonal subspaces of $L_2(P_3)$.
We assume that the parametrization is smooth in the following sense.
There is a linear space $K$, the \emph{tangent space} of $\Delta$,
and a bounded linear operator $D=(D_Q,D_R):K\to V\times W$,
and for each $k\in K$ there is a sequence $\delta_{nk}$ in $\Delta$
such that $Q_{nk}=Q_{\delta_{nk}}$ is Hellinger differentiable at $Q$
with derivative $D_Q k\in V$,
\[
P_1\Big[\int\Big(dQ_{nk}^{1/2}-dQ^{1/2}
-\frac{1}{2}n^{-1/2}D_Q k\,dQ^{1/2}\Big)^2\Big] \to 0,
\]
and $R_{nk}=R_{\delta_{nk}}$ is Hellinger differentiable at $R$
with derivative $D_R k\in W$,
\[
P_2\Big[\int\Big(dR_{nk}^{1/2}-dR^{1/2}
-\frac{1}{2}n^{-1/2}D_R k\,dR^{1/2}\Big)^2\Big] \to 0.
\]
Now write $M_n$ for the distribution of $Z_t$, $0\leq t\leq n$,
if $Q$ and $R$ are in effect, and $M_{nk}$ if $Q_{nk}$ and $R_{nk}$ are.
By Taylor expansion and \Ref{lln} and \Ref{clt},
we obtain \emph{local asymptotic normality}:
\bel{lan1}
\log\frac{dM_{nk}}{dM_n} = n^{-1/2}\sj
\big(D_Q k(X_{j-1},X_j)+D_R k(X_{j-1},X_j,U_j)\big)
- m^{-1}(P_2[D_Q^2 k] + P_3[D_R^2 k]) + \op
\ee
and
\bel{lan2}
n^{-1/2}\sj\big(D_Q k(X_{j-1},X_j)+D_R k(X_{j-1},X_j,U_j)\big)
\Rightarrow m^{-1/2}(P_2[D_Q^2 k] + P_3[D_R^2 k])^{1/2} Y.
\ee
For Markov chains, different proofs are in Penev \cite{Pe91},
Bickel \cite{Bi93} and Greenwood and Wefelmeyer \cite{GW95};
see also Bickel and Kwon \cite{BK01}. For Markov step processes
see  H\"opfner, Jacod and Ladelli \cite{HJL90} and
H\"opfner \cite{Ho93a, Ho93b}. A proof for nonparametric semi-Markov models
is in Greenwood and Wefelmeyer \cite{GW96}.

We want to estimate a $d$-dimensional functional $\f:\Delta\to\R^d$
of the parameter $\delta$. We call $\f$ \emph{differentiable} at $\delta$
with \emph{gradient} $(v_\f,w_\f)$ if $v_\f\in V^d$, $w_\f\in W^d$, and
\bel{grad}
n^{1/2}(\f(\delta_{nk})-\f(\delta))
\to m^{-1}(P_2[v_\f D_Q k]+P_3[w_\f D_R k]), \quad k\in K.
\ee
The \emph{canonical gradient} $(v_\f^*,w_\f^*)$ of $\f$ is the componentwise
projection of $(v_\f,w_\f)$ onto the closure of $(DK)^d$ in $(L_2(P_3))^d$.
If $DK$ is closed in $L_2(P_3)$, we can write
$(v_\f^*,w_\f^*)=(D_Q k_\f,D_R k_\f)$  for some $k_\f\in K$.
This will be the case in Sections 3--5.

An estimator $\hf$ is called \emph{regular} for $\f$ at $\delta$
with \emph{limit} $L$ if $L$ is a $d$-dimensional random vector such that
\[
n^{1/2}(\hf-\f(\delta_{nk})) \Rightarrow L \quad \mbox{under } M_{nk},
\quad k\in K.
\]

The convolution theorem says that
\[
L = A + m^{-1/2}(P_2[v_\f^* v_\f^{*\top}]
+ P_3[w_\f^* w_\f^{*\top}])^{1/2}Y_d,
\]
with $Y_d$ a $d$-dimensional standard normal random vector,
and $A$ a $d$-dimensional random vector independent of $Y_d$.
This justifies calling $\hf$ \emph{efficient} for $\f$ at $\delta$ if
$n^{1/2}(\hf-\f(\delta))$ is asymptotically normal under $M_n$
with covariance matrix
$m^{-1} (P_2[v_\f^* v_\f^{*\top}] + P_3[w_\f^* w_\f^{*\top}])$.

An estimator $\hf$ is called \emph{asymptotically linear} for $\f$ at $\delta$
with \emph{influence function} $(a,b)$ if $a\in V^d$, $b\in W^d$, and
\[
n^{1/2}(\hf-\f(\delta)) = n^{-1/2}\sj
\big(a(X_{j-1},X_j)+b(X_{j-1},X_j,U_j)\big) + \op.
\]
We have the following characterization. An estimator $\hf$ is regular
and efficient for $\f$ at $\delta$ if and only if it is asymptotically linear
with influence function equal to the canonical gradient,
\[
n^{1/2}(\hf-\f(\delta)) = n^{-1/2}\sj
\big(v_\f^*(X_{j-1},X_j)+w_\f^*(X_{j-1},X_j,U_j)\big) + \op.
\]
For proofs of the convolution theorem and the characterization we refer to
Bickel, Klaassen, Ritov and Wellner \cite {BKRW98}.

To prove asymptotic linearity of estimators in misspecified models,
we need the following \emph{martingale approximation}.
Set $L_{2,0}(P_2) = \{f\in L_2(P_2):P_2[f]=0\}$.
The \emph{potential} $G$ of the embedded Markov chain is defined by
\[
Gf = \sum_{i=0}^\infty Q^i f, \quad f\in L_{2,0}(P_2).
\]
For $f\in L_2(P_2)$ set
\[
Af(x,y) = G(f-P_2[f])(y) - QG(f-P_2[f])(x)
= \sum_{i=0}^\infty (Q^i f(y) - Q^{i+1}f(x)).
\]
Then $QAf=0$ and
\[
P_2[(Af)^2] = P_2[f^2] - (P_2[f])^2
+ 2\sum_{i=1}^\infty P_2[(f-P_2[f])Q^i f].
\]

Let $f\in L_2(P_3)$ and set $f_0=f-Rf$.
Then we obtain the stochastic expansion
\bel{mart}
n^{-1/2}\sj (f(X_{j-1},X_j,U_j)-P_3[f])
= n^{-1/2}\sj \big(ARf(X_{j-1},X_j)+f_0(X_{j-1},X_j,U_j)\big) + \op.
\ee
Note that $QAR f=0$ and $Sf_0=0$. Hence $ARf(X_{j-1},X_j)$
and $f_0(X_{j-1},X_j,U_j)$ are orthogonal martingale increments.
For discrete-time processes,
the martingale approximation \Ref{mart} is due to Gordin \cite{Go69}
and Gordin and Lif\v{s}ic \cite{GL78}. It was discovered independently
by Maigret \cite{Ma78}, D\"urr and Goldstein \cite{DG86}
and Greenwood and Wefelmeyer \cite{GW95}. See also Section 17.4
in the monograph of Meyn and Tweedie \cite{MT93}.
The martingale approximation \Ref{mart}
and the martingale central limit theorem \Ref{clt} imply that
\[
n^{-1/2}\sj (f(X_{j-1},X_j,U_j)-P_3[f])
\Rightarrow m^{-1/2}(P_2[(ARf)^2] + P_3[(f-Rf)^2])^{1/2}Y.
\]

To calculate canonical gradients of functionals in misspecified models,
we need the following \emph{perturbation expansion},
due to Kartashov \cite{Ka85a, Ka85b, Ka96},
\bel{per}
n^{1/2}(P_{2nk}[f]-P_2[f]) \to P_2[D_Q k \cdot Af], \quad k\in K.
\ee
Here $P_{2nk}$ denotes the distribution of $(X_{j-1},X_j)$ if $Q_{nk}$
is in effect. This pathwise version of the perturbation expansion
suffices for our purposes. Greenwood and Wefelmeyer \cite{GW95} show
that it follows also from the martingale approximation \Ref{mart}.

\section{Model Q}

In Model Q we assume a parametric model $q_\vt$, $\vt\in \Theta\subset\R^d$,
for the $\mu$-density of the transition distribution
of the embedded Markov chain, and consider the conditional inter-arrival time
distribution as unknown. Suppose the model is misspecified,
and the true transition distribution is $Q$. Then the KL functional $K_Q(P_2)$
maximizes $P_2[\log q_\vt]$, and the partial maximum likelihood estimator
$\hvt_Q$ maximizes $\P_2[\log q_\vt]$.
Write
\[
\chi_\vt(x,y)=\partial_\vt\log q_\vt(x,y)
\]
for the $d$-dimensional vector of partial derivatives of $\log q_\vt(x,y)$.
Then $K_Q(P_2)$ solves $P_2[\chi_\vt]=0$,
and $\hvt_Q$ solves $\P_2[\chi_\vt]=0$.
Heuristically, by Taylor expansion,
\ben
0 &=& \P_2[\chi_{\hvt_Q}] = \avj\chi_{\hvt_Q}(X_{j-1},X_j) \nonumber \\
\label{chi}
&=& \avj \chi_{K_Q(P_2)}(X_{j-1},X_j)
+ \avj \dc_{K_Q(P_2)}(X_{j-1},X_j)(\hvt_Q-K_Q(P_2)) + \opn.
\non
Here $\dc_\vt(x,y)$ is the $d\times d$ matrix of partial derivatives of
$\chi_\vt(x,y)$. With \Ref{N} and \Ref{lln} we obtain
\bel{i1}
n^{1/2}(\hvt_Q-K_Q(P_2))
= - m(P_2[\dc_{K_Q(P_2)}])^{-1}n^{-1/2}\sj\chi_{K_Q(P_2)}(X_{j-1},X_j) + \op.
\ee

If Model Q is correctly specified and $Q=Q_\vt$, then $K_Q(P_2)=\vt$.
We also have the following relations, which are well-known in the i.i.d.\ case,
\[
0 = \partial_\vt Q_\vt(\cdot,E) = Q_\vt\chi_\vt, \qquad
0 = \partial_\vt Q_\vt\chi_\vt = Q_\vt\chi_\vt \chi_\vt^\top + Q_\vt\dc_\vt.
\]
In particular, the partial Fisher information matrix for Model Q is
$I_\vt=-P_2[\dc_\vt]=P_2[\chi_\vt \chi_\vt^\top]$.
Hence, for the correctly specified model, the partial maximum
likelihood estimator $\hvt_Q$ has the stochastic expansion
\[
n^{1/2}(\hvt_Q-\vt) = mI_\vt^{-1}n^{-1/2}\sj\chi_\vt(X_{j-1},X_j) + \op.
\]
This means that $\hvt_Q$ is asymptotically linear with influence function
$mI_\vt^{-1}(\chi_\vt,0)$, and $n^{1/2}(\hvt_Q-\vt)$ is asymptotically normal
with covariance matrix $mI_\vt^{-1}$.

If the model is misspecified, then $\chi_{K_Q(P_2)}$ is not in $V^d$.
We apply the martingale approximation \Ref{mart} to \Ref{i1}
and see that $\hvt_Q$ is asymptotically linear with influence function
$-m(P_2[\dc_{K_Q(P_2)}])^{-1}(A\chi_{K_Q(P_2)},0)$.
Hence $n^{1/2}(\hvt_Q-K_Q(P_2))$ is asymptotically normal
with covariance matrix
\[
m(P_2[\dc_{K_Q(P_2)}])^{-1} P_2[A\chi_{K_Q(P_2)}A^\top\chi_{K_Q(P_2)}]
(P_2[\dc_{K_Q(P_2)}])^{-1}.
\]

Let us now prove efficiency of $\hvt_Q$,
first for the correctly specified model.
For $c\in\R^d$ set $\vt_{nc}=\vt+n^{-1/2}c$. Assume that $q_{nc}=q_{\vt_{nc}}$
is Hellinger differentiable at $\vt$,
\bel{hellq}
\int\int\Big(q_{nc}^{1/2}(x,y)-q_\vt^{1/2}(x,y)
-\frac{1}{2}n^{-1/2}c^\top\chi_\vt(x,y)q_\vt^{1/2}(x,y)\Big)^2
\mu(dy)P_1(dx) \to 0.
\ee
Let $\cR$ denote the set of all conditional inter-arrival distributions.
For $w\in W$ choose a sequence $R_{nw}$ in $\cR$
that is Hellinger differentiable at $R$,
\bel{R}
P_2\Big[\int\Big(dR_{nw}^{1/2}-dR^{1/2}
-\frac{1}{2}n^{-1/2}w\,dR^{1/2}\Big)^2\Big] \to 0.
\ee
Then the assumptions of Section 2 hold with $\Delta=\Theta\times\cR$,
$K=\R^d\times W$, $D_Q(c,w)=c^\top\chi_\vt$, $D_R(c,w)=w$.
The functional to be estimated is $\f(\vt,R)=\vt$.
By orthogonality of $V$ and $W$, its canonical gradient is
obtained from \Ref{grad} as $(c_\vt^\top\chi_\vt,0)$
with $d\times d$ matrix $c_\vt$ determined by
\[
c = m^{-1}c_\vt^\top P_2[\chi_\vt \chi_\vt^\top] c
= m^{-1}c_\vt^\top I_\vt c, \quad c\in\R,
\]
i.e.\ $c_\vt=mI_\vt^{-1}$. Hence the canonical gradient of $\vt$ is
$mI_\vt^{-1}(\chi_\vt,0)$ and equals the influence function of $\hvt$,
which is therefore efficient for the correctly specified model.

Suppose now that the model is misspecified, and let $\cQ$ be the set
of all transition distributions of the embedded Markov chain.
Let $Q$ denote the true transition distribution.
For $v\in V$ choose a sequence $Q_{nv}$ in $\cQ$ that is
Hellinger differentiable at $Q$,
\bel{Q}
P_1\Big[\int\Big(dQ_{nv}^{1/2}-dQ^{1/2}
-\frac{1}{2}n^{-1/2}v\,dQ^{1/2}\Big)^2\Big] \to 0.
\ee
Then the assumptions of Section 2 hold with $\Delta=\cQ\times\cR$,
$K=V\times W$, $D_Q(v,w)=v$, $D_R(v,w)=w$.
The functional to be estimated is $\f(Q,R)=K_Q(P_2)$. Heuristically,
\[
0 = P_{2nv}[\chi_{K_Q(P_{2nv})}] = P_{2nv}[\chi_{K_Q(P_2)}]
+ P_{2nv}[\dc_{K_Q(P_2)}](K_Q(P_{2nv})-K_Q(P_2)) + \opn.
\]
With $P_{2nv}[\dc_{K_Q(P_2)}]\to P_2[\dc_{K_Q(P_2)}]$ we obtain
\[
K_Q(P_{2nv})-K_Q(P_2) = -(P_2[\dc_{K_Q(P_2)}])^{-1}P_{2nv}[\chi_{K_Q(P_2)}]
+ \opn.
\]
The perturbation expansion \Ref{per} yields
\bel{perq}
n^{1/2}P_{2nv}[\chi_{K_Q(P_2)}]
= n^{1/2}(P_{2nv}-P_2)[\chi_{K_Q(P_2)}] \to P_2[vA\chi_{K_Q(P_2)}].
\ee
Hence
\[
n^{1/2}(K_Q(P_{2nv})-K_Q(P_2))
\to - (P_2[\dc_{K_Q(P_2)}])^{-1} P_2[vA\chi_{K_Q(P_2)}], \quad v\in V,
\]
and the canonical gradient of $K_Q$ is obtained from \Ref{grad} as
$-m(P_2[\dc_{K_Q(P_2)}])^{-1}(A\chi_{K_Q(P_2)},0)$
and equals the influence function of $\hvt_Q$,
which is therefore efficient for the misspecified model.

\section{Model R}

Model R is completely analogous to Model Q, with interchanged roles
of the transition distribution $Q$ of the embedded Markov chain,
and the  conditional inter-arrival time distribution $R$. Specifically,
in Model R we assume a parametric model $r_\vt$, $\vt\in \Theta\subset\R^d$,
for the $\nu$-density of the conditional inter-arrival time,
and consider the transition distribution of the embedded Markov chain
as unknown. Suppose the model is misspecified,
and the true conditional inter-arrival time distribution is $R$.
Then the KL functional $K_R(P_3)$ maximizes $P_3[\log r_\vt]$,
and the partial maximum likelihood estimator $\hvt_Q$ maximizes
$\P_3[\log r_\vt]$. Write
\[
\vr_\vt(x,y,u)=\partial_\vt\log r_\vt(x,y,u)
\]
for the $d$-dimensional vector of partial derivatives of $\log r_\vt(x,y,u)$.
Then $K_R(P_3)$ solves $P_3[\vr_\vt]=0$, and $\hvt_R$ solves $\P_3[\vr_\vt]=0$.
Heuristically, by Taylor expansion,
\ben
0 &=& \P_3[\vr_{\hvt_R}] = \avj\vr_{\hvt_R}(X_{j-1},X_j,U_j) \nonumber \\
\label{vr}
&=& \avj \vr_{K_R(P_3)}(X_{j-1},X_j,U_j)
+ \avj \dr_{K_R(P_3)}(X_{j-1},X_j,U_j)(\hvt_R-K_R(P_3)) + \opn.
\non
Here $\dr_\vt(x,y,u)$ is the $d\times d$ matrix of partial derivatives of
$\vr_\vt(x,y,u)$. With \Ref{N} and \Ref{lln} we obtain
\bel{i2}
n^{1/2}(\hvt_R-K_R(P_3)) = - m(P_3[\dr_{K_R(P_3)}])^{-1}n^{-1/2}
\sj\vr_{K_R(P_3)}(X_{j-1},X_j,U_j) + \op.
\ee

If Model R is correctly specified and $R=R_\vt$, then $K_R(P_3)=\vt$.
We also have the following relations,
\[
0 = \partial_\vt R_\vt(\cdot,\cdot,\R) = R_\vt\vr_\vt, \qquad
0 = \partial_\vt R_\vt\vr_\vt = R_\vt\vr_\vt \vr_\vt^\top + R_\vt\dr_\vt.
\]
In particular, the partial Fisher information matrix for Model R is
$J_\vt=-P_3[\dr_\vt]=P_3[\vr_\vt \vr_\vt^\top]$.
Hence, for the correctly specified model, the partial maximum
likelihood estimator $\hvt_R$ has the stochastic expansion
\[
n^{1/2}(\hvt_R-\vt) = mJ_\vt^{-1}n^{-1/2}\sj\vr_\vt(X_{j-1},X_j,U_j) + \op.
\]
This means that $\hvt_R$ is asymptotically linear with influence function
$mJ_\vt^{-1}(0,\vr_\vt)$, and $n^{1/2}(\hvt_R-\vt)$ is asymptotically normal
with covariance matrix $mJ_\vt^{-1}$.

If the model is misspecified, then $\vr_{K_R(P_3)}$ is not in $W^d$.
We apply the martingale approximation \Ref{mart} to \Ref{i2}
and see that $\hvt_R$ is asymptotically linear with influence function
\[
-m(P_3[\dr_{K_R(P_3)}])^{-1}
(AR\vr_{K_R(P_3)},\vr_{K_R(P_3)}-R\vr_{K_R(P_3)}).
\]
Hence $n^{1/2}(\hvt_R-K_R(P_3))$ is asymptotically normal with covariance
matrix
\[
m(P_3[\dr_{K_R(P_3)}])^{-1}\Sigma_R(P_3[\dr_{K_R(P_3)}])^{-1},
\]
where
\[
\Sigma_R = P_2[AR\vr_{K_R(P_3)}A^\top R\vr_{K_R(P_3)}]
+P_3[(\vr_{K_R(P_3)}-R\vr_{K_R(P_3)})(\vr_{K_R(P_3)}-R\vr_{K_R(P_3)})^\top].
\]

Let us now prove efficiency of $\hvt_R$,
first for the correctly specified model.
For $c\in\R^d$ set $\vt_{nc}=\vt+n^{-1/2}c$. Assume that $r_{nc}=r_{\vt_{nc}}$
is Hellinger differentiable at $\vt$,
\bel{hellr}
\int\int\Big(r_{nc}^{1/2}(x,y,u)-r_\vt^{1/2}(x,y,u)
-\frac{1}{2}n^{-1/2}c^\top\vr_\vt(x,y,u)r_\vt^{1/2}(x,y,u)\Big)^2
\nu(du)P_2(d(x,y)) \to 0.
\ee
Let $\cQ$ denote the set of all transition distributions
of the embedded Markov chain.
For $v\in V$ choose a sequence $Q_{nv}$ in $\cQ$
that is Hellinger differentiable \Ref{Q} at $Q$.
Then the assumptions of Section 2 hold with $\Delta=\cQ\times\Theta$,
$K=V\times\R^d$, $D_Q(v,c)=v$, $D_R(v,c)=c^\top\vr_\vt$.
The functional to be estimated is $\f(Q,\vt)=\vt$.
By orthogonality of $V$ and $W$, its canonical gradient is
obtained from \Ref{grad} as $(0,c_\vt^\top\vr_\vt)$
with $d\times d$ matrix $c_\vt$ determined by
\[
c = m^{-1}c_\vt^\top J_\vt c, \quad c\in\R,
\]
i.e.\ $c_\vt=mJ_\vt^{-1}$. Hence the canonical gradient of $\vt$ is
$mJ_\vt^{-1}(0,\vr_\vt)$ and equals the influence function of $\hvt$,
which is therefore efficient for the correctly specified model.

Suppose now that the model is misspecified, and let $\cR$ be the set
of all transition distributions of the embedded Markov chain.
Let $R$ denote the true transition distribution.
For $w\in W$ choose a sequence $R_{nw}$ in $\cR$ that is
Hellinger differentiable \Ref{R} at $R$.
Then the assumptions of Section 2 hold with $\Delta=\cQ\times\cR$,
$K=V\times W$, $D_Q(v,w)=v$, $D_R(v,w)=w$.
The functional to be estimated is $\f(Q,R)=K_R(P_3)$. Heuristically,
\[
0 = P_{3nvw}[\vr_{K_R(P_{3nvw})}] = P_{3nvw}[\vr_{K_R(P_3)}]
+ P_{3nvw}[\dr_{K_R(P_3)}](K_R(P_{3nvw})-K_R(P_3)) + \opn.
\]
With $P_{3nvw}[\dr_{K_R(P_3)}]\to P_3[\dr_{K_R(P_3)}]$ we obtain
\[
K_R(P_{3nvw})-K_R(P_3) = -(P_3[\dr_{K_R(P_3)}])^{-1}P_{3nvw}[\vr_{K_R(P_3)}]
+ \opn.
\]
Write $P_{3nvw}=P_{2nv}\otimes R_{nw}$ and apply the perturbation
expansion \Ref{perq} to obtain
\be*
n^{1/2}(K_R(P_{3nvw})-K_R(P_3))
&\to& - (P_3[\dr_{K_R(P_3)}])^{-1}
\Big(P_2[vAR\vr_{K_R(P_3)}]+P_3[w\vr_{K_R(P_3)}]\Big) \\
&& = - (P_3[\dr_{K_R(P_3)}])^{-1} \Big(P_2[vAR\vr_{K_R(P_3)}]
+P_3[w(\vr_{K_R(P_3)}-R\vr_{K_R(P_3)})]\Big),
\no*
and the canonical gradient of $K_R$ is obtained from \Ref{grad} as
\[
-m(P_3[\dr_{K_R(P_3)}])^{-1}(AR\vr_{K_R(P_3)},
\vr_{K_R(P_3)}-R\vr_{K_R(P_3)})
\]
and equals the influence function of $\hvt_R$,
which is therefore efficient for the misspecified model.

\section{Model S}

While Models Q and R are semiparametric, Models S is parametric.
In Model S we assume parametric models $q_\vt$ and $r_\vt$,
$\vt\in \Theta\subset\R^d$,
for the $\mu$-density of the transition distribution
of the embedded Markov chain and for the $\nu$-density of the conditional
inter-arrival time. We have $s_\vt(x,y,u)=q_\vt(x,y)r_\vt(x,y,u)$.
Hence the KL functional $K_S(P_3)$ maximizes
$P_3[\log s_\vt]=P_2[\log q_\vt]+P_3[\log r_\vt]$,
and the partial maximum likelihood estimator $\hvt_S$ maximizes
$\P_3[\log s_\vt]=\P_2[\log q_\vt]+\P_3[\log r_\vt]$.
Write
\[
\sigma_\vt(x,y,u)=\partial_\vt\log s_\vt(x,y,u)=\chi_\vt(x,y)+\vr_\vt(x,y,u)
\]
for the $d$-dimensional vector of partial derivatives of $\log s_\vt(x,y,u)$.
Then $K_S(P_3)$ solves $P_3[\sigma_\vt]=P_2[\chi_\vt]+P_3[\vr_\vt]=0$,
and $\hvt_S$ solves $\P_3[\sigma_\vt]=\P_2[\chi_\vt]+\P_3[\vr_\vt]=0$.
Taylor expansions analogous to \Ref{chi} and \Ref{vr} imply
\be*
0 &=& \P_3[\sigma_{\hvt_S}] = \avj\sigma_{\hvt_S}(X_{j-1},X_j,U_j) \nonumber \\
&=& \avj \sigma_{K_S(P_3)}(X_{j-1},X_j,U_j)
+ \avj \ds_{K_S(P_3)}(X_{j-1},X_j,U_j)(\hvt_S-K_S(P_3)) + \opn,
\no*
where $\ds_\vt(x,y,u)=\dc_\vt(x,y)+\dr_\vt(x,y,u)$
is the $d\times d$ matrix of partial derivatives of $\sigma_\vt(x,y,u)$.
We obtain
\bel{i3}
n^{1/2}(\hvt_S-K_S(P_3))
= - m(P_3[\ds_{K_S(P_3)}])^{-1}
n^{-1/2}\sj\sigma_{K_S(P_3)}(X_{j-1},X_j,U_j) + \op.
\ee

If Model S is correctly specified with $Q=Q_\vt$ and $R=R_\vt$,
then $K_S(P_3)=\vt$. From Sections 3 and 4 we obtain the Fisher information
matrix for Model S as $I_\vt+J_\vt$.
Hence, for the correctly specified model, the maximum
likelihood estimator $\hvt_S$ has the stochastic expansion
\[
n^{1/2}(\hvt_S-\vt) = m(I_\vt+J_\vt)^{-1}
n^{-1/2}\sj\sigma_\vt(X_{j-1},X_j,U_j) + \op.
\]
This means that $\hvt_S$ is asymptotically linear with influence function
$m(I_\vt+J_\vt)^{-1}(\chi_\vt,\vr_\vt)$, and $n^{1/2}(\hvt_S-\vt)$
is asymptotically normal with covariance matrix $m(I_\vt+J_\vt)^{-1}$.

If the model is misspecified, then $\chi_{K_S(P_3)}$ is not in $V^d$
and $\vr_{K_S(P_3)}$ is not in $W^d$.
We apply the martingale approximation \Ref{mart} to \Ref{i3}
and see that $\hvt_S$ is asymptotically linear with influence function
\[
-m(P_3[\ds_{K_S(P_3)}])^{-1}
(A\chi_{K_S(P_3)}+AR\vr_{K_S(P_3)},\vr_{K_S(P_3)}-R\vr_{K_S(P_3)}).
\]
Hence $n^{1/2}(\hvt_S-K_S(P_3))$ is asymptotically normal
with covariance matrix
\[
m(P_3[\ds_{K_S(P_3)}])^{-1}\Sigma_S(P_3[\ds_{K_S(P_3)}])^{-1},
\]
where
\[
\Sigma_S = P_2[A(\chi_{K_S(P_3)}+R\vr_{K_S(P_3)})
A^\top(\chi_{K_S(P_3)}+R\vr_{K_S(P_3)})]
+ P_3[(\vr_{K_S(P_3)}-R\vr_{K_S(P_3)})(\vr_{K_S(P_3)}-R\vr_{K_S(P_3)})^\top].
\]

Let us now prove efficiency of $\hvt_S$,
first for the correctly specified model.
For $c\in\R^d$ set $\vt_{nc}=\vt+n^{-1/2}c$. Assume that $q_{nc}=q_{\vt_{nc}}$
is Hellinger differentiable \Ref{hellq} at $\vt$, and $r_{nc}=r_{\vt_{nc}}$
is Hellinger differentiable \Ref{hellr} at $\vt$.
Then the assumptions of Section 2 hold with $\Delta=\Theta$, $K=\R^d$,
$D_Q c=c^\top\chi_\vt$, $D_R c=c^\top\vr_\vt$.
The functional to be estimated is $\f(\vt)=\vt$.
The canonical gradient is obtained from \Ref{grad} as
$m(I_\vt+J_\vt)^{-1}(\chi_\vt,\vr_\vt)$. It equals the influence function
of $\hvt_S$, which is therefore efficient in the correctly specified model.

Suppose now that the model is misspecified. Let $\cQ$ be the set
of all transition distributions of the embedded Markov chain,
and let $\cR$ be the set of all transition distributions of the embedded
Markov chain. For $v\in V$ choose a sequence $Q_{nv}$ in $\cQ$
that is Hellinger differentiable \Ref{Q} at $Q$.
For $w\in W$ choose a sequence $R_{nw}$ in $\cR$ that is
Hellinger differentiable \Ref{R} at $R$.
Then the assumptions of Section 2 hold with $\Delta=\cQ\times\cR$,
$K=V\times W$, $D_Q(v,w)=v$, $D_R(v,w)=w$.
The functional to be estimated is $\f(Q,R)=K_S(P_3)$.
Similarly as in Section 4,
\be*
& 0 = P_{3nvw}[\vr_{K_S(P_{3nvw})}] = P_{3nvw}[\vr_{K_S(P_3)}]
+ P_{3nvw}[\dr_{K_S(P_3)}](K_S(P_{3nvw})-K_S(P_3)) + \opn, & \\
& K_S(P_{3nvw})-K_S(P_3)
= -(P_3[\ds_{K_S(P_3)}])^{-1}P_{3nvw}[\sigma_{K_S(P_3)}]
+ \opn,
\no*
and therefore
\[
n^{1/2}(K_S(P_{3nvw})-K_R(P_3))
\to - (P_3[\dr_{K_S(P_3)}])^{-1} \Big(P_2[v(A\chi_{K_S(P_3)}+AR\vr_{K_S(P_3)})]
+P_3[w(\vr_{K_S(P_3)}-R\vr_{K_S(P_3)})]\Big).
\]
Hence by \Ref{grad} the canonical gradient of $K_S$ is obtained as
\[
-m(P_3[\ds_{K_S(P_3)}])^{-1}(A\chi_{K_S(P_3)}+AR\vr_{K_S(P_3)},
\vr_{K_S(P_3)}-R\vr_{K_S(P_3)})
\]
and equals the influence function of $\hvt_S$,
which is therefore efficient for the misspecified model.

\section{Remarks}

In this section we comment on examples and possible extensions
of our results.

\textbf{1.} If the distribution of the inter-arrival times charges only 1,
so that $R(x,y,du)=\delta_1(du)$, then the semi-Markov process reduces
to a Markov chain with transition distribution $Q$, and for Model Q
we recover the results of Greenwood and Wefelmeyer \cite{GW97}.

\textbf{2.} Our results carry over to observations $(X_0,T_0),\dots,(X_n,T_n)$
of the embedded Markov renewal process. Just replace $N$ by $n$.
In particular, instead of the central limit theorem \Ref{clt}
with random summation index $N$, use
\[
n^{-1/2}\sum_{j=1}^n f(X_{j-1},X_j,U_j) \Rightarrow (P_3[f^2])^{1/2} Y,
\]
and replace $m$ by 1 everywhere.

\bigskip\noindent
In some examples we can describe the KL functional more explicitly.

\textbf{3.} Suppose the embedded Markov chain is a linear autoregressive model
of order 1, i.e.\ $X_j=\vt X_{j-1}+\ve_j$, where $\vt\in\R$
and the innovations $\ve_j$ are i.i.d.\ with mean 0, finite variance,
and known density $f$. Then Model Q holds with $Q(x,dy)=f(y-\vt x)dy$,
and $\chi_\vt(x,y)=x\ell(y-\vt x)$ with $\ell=-f'/f$.
Hence the KL functional solves $E[X_0\ell(X_1-\vt X_0)]=0$.
If $f$ is the density of $\tau Y$ for some $\tau>0$,
then $\ell(x)=\tau^{-2}x$ and
$E[X_0\ell(X_1-\vt X_0)]=\tau^{-2}(E[X_0 X_1]-\vt E[X_0^2])$.
Hence the KL functional is $K_Q(P_2)=E[X_0 X_1]/E[X_0^2]$,
and the partial maximum likelihood estimator for $\vt$ is
the least squares estimator
\[
\hvt_Q = K_Q(\P_2) = \frac{\sj X_{j-1}X_j}{\sj X_{j-1}^2},
\]
a ratio of two empirical estimators.

\textbf{4.} Suppose the inter-arrival time $U_j$ given $X_{j-1}=x$ and $X_j=y$
is exponentially distributed with mean $1/\lambda(x)$ not depending on $y$,
\[
R(x,y,du) = \lambda(x)\exp(-u\lambda(x))du.
\]
Then the semi-Markov process is a Markov step process.
If the mean is constant, $\lambda(x)=\vt$, $\vt>0$, then Model R holds
with $R_\vt(x,y,du)=\vt\exp(\vt u)$, and $\vr_\vt(x,y,u)=\vt^{-1}-u$.
Hence the KL functional solves $E[\vr_\vt(X_0,X_1,U_1)]=\vt^{-1}-E[U_1]=0$,
and we obtain $K_R(P_3)=1/E[U_1]$. The partial maximum likelihood estimator
for $\vt$ is
\[
\hvt_R = 1\Big/\avj U_j,
\]
a function of an empirical estimator. Efficiency of empirical estimators
in Markov step processes is studied in Greenwood and Wefelmeyer \cite{GW94}.

\bigskip\noindent
The models Q, R and S are described in terms of the \emph{conditional}
distributions $Q(x,dy)$ and $R(x,y,du)$. It is occasionally reasonable
to model instead the \emph{marginal} distributions $P_1$, $P_2$ or $P_3$.
Results for these three models differ considerably among each other
and from Models Q, R and S.

\textbf{5.} Suppose we have a parametric model for the $\mu$-density
$p_{1\vt}$ of $P_1$. The \emph{marginal maximum likelihood estimator} $\hvt_1$
maximizes
\[
\P_1[\log p_{1\vt}] = \avj\log p_{1\vt}(X_{j-1}).
\]
It estimates the \emph{KL functional} $K(P_1)$, the parameter
that maximizes $P_1[\log p_{1\vt}]$. Note that
the marginal maximum likelihood estimator is an empirical version
of the KL functional, $\hvt_1=K(\P_1)$.

However, $\hvt_1$ is \emph{not} efficient for $\vt$ when the marginal model
is correctly specified.
The reason is that the specification $p_{1\vt}$ of the marginal density
implies a constraint on the conditional distribution $Q$ of the embedded
Markov chain, but the marginal maximum likelihood estimator does not use
this information. An efficient estimator for $\vt$ is difficult to construct.
See Kessler, Schick and Wefelmeyer \cite{KSW01} for an efficient estimator
of $\vt$ in a Markov chain model with a (correctly specified) parametric
model for the (one-dimensional) marginal density.
On the other hand, $\hvt_1$ is efficient for $K(P_1)$
in a nonparametric sense when the marginal model is misspecified.

We note that, in this respect, semi-Markov processes and Markov chains
are different from the i.i.d.\ case. Suppose we have i.i.d.\ observations
$(X_j,Y_j)$ with joint distribution $p_{1\vt}(x)dx\,Q(x,dy)$,
where $Q$ is unknown. Then $Q$ is not constrained by the marginal model
$p_{1\vt}$, and the marginal maximum likelihood estimator is efficient
for $\vt$ if the marginal model is correctly specified,
and also efficient for $K(P_1)$ if the marginal model is misspecified.

\textbf{6.} Suppose we have a parametric model for the $\mu^2$-density
$p_{2\vt}$ of $P_2$. The \emph{marginal maximum likelihood estimator} $\hvt_2$
maximizes
\[
\P_2[\log p_{2\vt}] = \avj\log p_{2\vt}(X_{j-1},X_j).
\]
It estimates the \emph{KL functional} $K(P_2)$, the parameter
that maximizes $P_2[\log p_{2\vt}]$, and $\hvt_2=K(\P_2)$.
The perturbation expansion \Ref{per} suggests that maximizing
$\P_2[\log p_{2\vt}]$ is asymptotically equivalent to solving
$\P_2[A\chi_\vt]=0$, and the martingale approximation
\Ref{mart} suggests that this is asymptotically equivalent
to solving $\P_2[\chi_\vt]=0$.
Hence the marginal maximum likelihood estimator $\hvt_2$
is asymptotically equivalent to the conditional maximum likelihood
estimator $\hvt_Q$ and therefore efficient in the correctly specified
model $p_{2\vt}$. The reason is that
$p_{2\vt}(x,y)=p_{1\vt}(x)q_\vt(x,y)$, and $q_\vt(x,y)$ determines
$p_{1\vt}$, which therefore does not contain additional information
about $\vt$.

This is again different from the i.i.d.\ case.
Suppose we have i.i.d.\ observations $(X_j,Y_j)$ with joint density
$p_{1\vt}(x)q_\vt(x,y)$. Then $p_{1\vt}$ contains, in general,
additional information about $\vt$.

\textbf{7.} Suppose we have a parametric model
for the $\mu^2\otimes\nu$-density $p_{3\vt}$ of $P_3$.
The \emph{marginal maximum likelihood estimator} $\hvt_3$
maximizes
\[
\P_3[\log p_{3\vt}] = \avj\log p_{3\vt}(X_{j-1},X_j,U_j).
\]
It estimates the \emph{KL functional} $K(P_3)$, the parameter
that maximizes $P_3[\log p_{3\vt}]$, and $\hvt_3=K(\P_3)$.
We can write $p_{3\vt}(x,y,u)=p_{2\vt}(x,y)r_\vt(x,y,u)$.
Now $r_\vt(x,y,u)$ carries additional information about $\vt$,
similarly as in the i.i.d.\ case.

\bigskip
\textbf{8.} Remark 5 tells us in particular the following,
rather obvious, fact.
If a parametric estimator is efficient in a nonparametric sense,
then the reason is not that it is efficient in a parametric model.
Rather, an estimator usually is nonparametrically efficient
because it is a smooth function of an empirical estimator.
We can illustrate this also with Model S. Suppose we have
parametric models $q_\vt$ and $r_\vt$ for the densities of $Q$ and $R$.
Let $\hvt_Q=K_Q(\P_2)$ be the conditional maximum likelihood estimator
based on the model $q_\vt$ alone. In general, $\hvt_Q$ will not be efficient
for $\vt$ when model S is correctly specified,
because $\hvt_Q$ does not use the information about $\vt$
in the model $r_\vt$.
But if both $q_\vt$ and $r_\vt$ are misspecified,
$\hvt_Q$ will be nonparametrically efficient for $K_Q(P_2)$,
which is the KL functional for Model Q but not for Model S.

\end{document}